\documentclass{amsart}
\usepackage{wasowicz}

\renewcommand{\C}{\mathcal{C}^2(S)}
\newcommand{\barp}{\overline{p}}
\newcommand{\df}{\|d^2f\|_{\infty}}
\newcommand{\dgpf}{\bigl\|d^2(g+f)\bigr\|_{\infty}}
\newcommand{\dgmf}{\bigl\|d^2(g-f)\bigr\|_{\infty}}
\renewcommand{\le}{\leqslant}
\renewcommand{\ge}{\geqslant}

\info{Published: Math. Inequal. Appl. 11/4 (2008), 693--700}

\begin{document}

\title{Hermite--Hadamard--type inequalities in the approximate integration}
\author{\SW}
\address{\SWaddr}
\email{\SWmail}
\date{October 31, 2007}
\keywords{%
 Approximate integration,
 Hermite--Hadamard inequality,
 convex functions,
 cubatures,
 norm of a quadratic form}
\subjclass[2000]{Primary: 26D15, 41A80, 65D32; Secondary: 15A63, 26B25, 39B62}
\begin{abstract}
 We give a slight extension of the Hermite--Hadamard inequality on simplices
 and we use it to establish error bounds of the operators connected with the
 approximate integration.
\end{abstract}
\maketitle

\section{Introduction}

In the numerical analysis there are many classical quadrature rules. For each
of them the error terms are well known. In the series of recent papers
\cite{Was07JIPAM,Was06JIPAM,Was08JIPAM2} the author established error bounds of quadrature operators
under regularity assumptions weaker from the classical ones. The method was
based on the use of Hermite--Hadamard--type inequalities for convex functions
of higher order.

In this paper we develop an analogous idea for convex functions of several
variables defined on simplices in $\R^n$. It is possible because the
generalizations of the Hermite--Hadamard inequality are known for multiple
integrals, e.g. by Choquet theory. Recently Bessenyei proved in~\cite{Bes08} such
an inequality on simplices using an elementary approach. We slightly extend
this result and we use this extension to give error bounds of operators
connected with the approximate computation of multiple integrals (i.e. the
cubatures).

Let $p_0,\dots,p_n\in\R^n$ be affine independent and let $S$ be a simplex with
vertices \mbox{$p_0,\dots,p_n$.} Denote by $\barp$ its barycenter, i.e.
\[
 \barp:=\frac{1}{n+1}\sum_{i=0}^np_i
\]
and by $\vol(S)$ its volume.

In our proofs we use some results contained in~\cite{Bes08}. One of them is the
following
\begin{lem}\label{Lm_aff}
 If $A:\R^n\to\R$ is an affine function then
 \[
  A(\barp)=\frac{1}{\vol(S)}\int_SA(x)dx=\frac{1}{n+1}\sum_{i=0}^nA(p_i).
 \]
\end{lem}

The properties of convex functions needed in this paper can be found in many
textbooks and monographs, for instance in \cite{Roc70} (cf. also~\cite{Bes08}).

Recall that a linear functional $\T$ defined on a linear space $X$ of (not
necessarily all) functions mapping some nonempty set into $\R$ is called
\emph{positive} if
\[
 f\le g\implies \T(f)\le\T(g)
\]
for any $f,g\in X$.

\section{Hermite--Hadamard--type inequalities}

We start with a slight generalization of the Hermite--Hadamard inequality.

\begin{thm}\label{Th_Had}
 Let $\T$ be positive linear functional defined (at least) on a~linear subspace
 of all functions mapping $S$ into $\R$ generated by a cone of convex
 functions. Assume that
 \[
  \T(A)=\frac{1}{\vol(S)}\int_SA(x)dx
 \]
 for any affine function $A:S\to\R$. Then the inequality
 \[
  f(\barp)\le\T(f)\le\frac{1}{n+1}\sum_{i=0}^nf(p_i)
 \]
 holds for any convex function $f:S\to\R$.
\end{thm}

\begin{proof}
 By convexity the subdifferential $\partial f(\barp)$ is non--empty,
 whence there exists an affine function $A:S\to\R$ supporting $f$ at the
 point $\barp$ (i.e. $A(\barp)=f(\barp)$ and $A\le f$ on $S$). Using
 Lemma~\ref{Lm_aff} and the properties of $\T$ we get
 \[
  f(\barp)=A(\barp)=\frac{1}{\vol(S)}\int_SA(x)dx=\T(A)\le\T(f).
 \]
 Since $p_0,\dots,p_n$ are affine independent, there exists (exactly one)
 affine function $B:S\to\R$ such that $B(p_i)=f(p_i)$, $i=0,\dots,n$. A simplex
 $S$ is the convex hull of the set $\{p_0,\dots,p_n\}$. Then by convexity
 $f\le B$ on $S$. Using once again Lemma~\ref{Lm_aff} and the properties of $\T$
 we arrive at
 \[
  \T(f)\le\T(B)=\frac{1}{\vol(S)}\int_SB(x)dx=\frac{1}{n+1}\sum_{i=0}^nB(p_i)
  =\frac{1}{n+1}\sum_{i=0}^nf(p_i).
 \]
\end{proof}

As an immediate consequence we obtain the Hermite--Hadamard inequality.
\begin{cor}
 If $f:S\to\R$ is convex then
 \[
  f(\barp)\le\frac{1}{\vol(S)}\int_Sf(x)dx\le\frac{1}{n+1}\sum_{i=0}^nf(p_i).
 \]
\end{cor}

Another important consequence of Theorem~\ref{Th_Had} concerns the operators of
the form of a convex combination of the values of $f$ at some appropriately
chosen points of $S$.

\begin{cor}
 Let $m\in\N$, $\lambda_1,\dots,\lambda_m\ge 0$ with $\sum_{i=1}^m\lambda_i=1$
 and $\xi_1,\dots,\xi_m\in S$ with $\sum_{i=1}^m\lambda_i\xi_i=\barp$. If a
 function $f:S\to\R$ is convex then
 \[
  f(\barp)\le\sum_{i=1}^m\lambda_if(\xi_i)\le\frac{1}{n+1}\sum_{i=0}^nf(p_i).
 \]
\end{cor}
\begin{proof}
 Clearly
 \[\T(f):=\sum_{i=1}^m\lambda_if(\xi_i)\]
 is a positive linear functional. Let $A:S\to\R$ be an affine function. Then by
 Lemma~\ref{Lm_aff}
 \[
  \T(A)=\sum_{i=1}^m\lambda_iA(\xi_i)=A\Bigl(\sum_{i=1}^m\lambda_i\xi_i\Bigr)
  =A(\barp)=\frac{1}{\vol(S)}\int_SA(x)dx
 \]
 and all the assumptions of Theorem~\ref{Th_Had} are fulfilled.
\end{proof}
This result is useful in the approximate computation of multiple integrals over
$S$. Together with further results of this paper it allows us to give error
bounds of the cubatures of the form of a conical combination of the values of
an integrand at some points of $S$. These error bounds depend on the second
order differential of the function involved.
\par\medskip
To conclude this section notice that in Theorem~\ref{Th_Had} not a form of an
operator $\T$ is essential but three its properties are important: $\T$ is
linear, positive and it coincides with the integral mean value in the class of
affine functions.

\section{A norm of a second order differential}

Let $\varphi:\R^n\to\R$ be a quadratic form, i.e.
\[
 \varphi(x)=\sum_{i,j=1}^na_{ij}x_ix_j
\]
for $x=(x_1,\dots,x_n)\in\R^n$. Let
\[
 \|\varphi\|:=\sup\bigl\{\bigl|\varphi(x)\bigr|:\|x\|=1\bigr\}.
\]
This is well defined norm on a linear space of all quadratic forms of $n$
variables with the following properties:
\begin{align*}
 \|\varphi\|&\le\sum_{i,j=1}^n|a_{ij}|,\\
 \bigl|\varphi(x)\bigr|&\le\|\varphi\|\cdot\|x\|^2,\quad x\in\R^n.
\end{align*}
The details may be easily checked (cf. e.g.~\cite{Sha96}).
\par\medskip
Now let $f\in\C$. Then the partial derivatives of the second order are bounded,
say, there exists a constant $M>0$ such that
\[
 \biggl|\frac{\partial^2f}{\partial x_i\partial x_j}(u)\biggr|\le\frac{M}{n^2}
\]
for any $u\in S$, $i,j=1,\dots,n$. Then for any $u\in S$
\[
 d^2f(u)(v,v)
 =\sum_{i,j=1}^n\frac{\partial^2f}{\partial x_i\partial x_j}(u)v_iv_j,
\]
whence
\[
 \bigl\|d^2f(u)\bigr\|\le\sum_{i,j=1}^n
 \biggl|\frac{\partial^2f}{\partial x_i\partial x_j}(u)\biggr|\le M.
\]
This allows us to define
\[
 \df:=\sup_{u\in S}\bigl\|d^2f(u)\bigr\|.
\]
Obviously
\begin{align}
 \bigl\|d^2(f+g)\bigr\|_{\infty}
 &\le\df+\|d^2g\|_{\infty},\label{norm1}\\
 \bigl\|d^2(\alpha f)\bigr\|_{\infty}&=|\alpha|\cdot\df,\label{norm2}\\
 \bigl|d^2f(u)(v,v)\bigr|&\le\bigl\|d^2f(u)\bigr\|\cdot\|v\|^2\le\df\|v\|^2
 \label{norm3}
\end{align}
for any $f,g\in\C$, $\alpha\in\R$, $u\in S$ and $v\in\R^n$.

\begin{lem}\label{Lm_norm}
 If $f\in\C$ and $g(x)=\frac{\df}{2}\|x\|^2$ then the functions
 $g+f$ and $g-f$ are convex. Moreover,
 \[
  \dgpf\le2\df\quad\text{and}\quad\dgmf\le2\df.
 \]
\end{lem}
\begin{proof}
 It is easy to compute $d^2g(u)(v,v)=\df\|v\|^2$. Then by~\eqref{norm3}
 \[
  \bigl|d^2f(u)(v,v)\bigr|\le\df\|v\|^2=d^2g(u)(v,v).
 \]
 Hence, for any $u\in S$, $v\in\R^n$,
 \[
  d^2(g+f)(u)(v,v)\ge 0\quad\text{and}\quad d^2(g-f)(u)(v,v)\ge 0,
 \]
 which means that for any $u\in S$, $d^2(g+f)(u)$ and $d^2(g-f)(u)$ are positive
 semi--definite. Then $g+f$ and $g-f$ are convex.
 \par\medskip\noindent
 If $\|v\|=1$ then $d^2g(u)(v,v)=\df$ for any $u\in S$, which implies
 $\bigl\|d^2g(u)\bigr\|=\df$, $u\in S$. Then $\|d^2g\|_{\infty}=\df$
 and by \eqref{norm1}, \eqref{norm2} we obtain the second assertion of the
 Lemma.
\end{proof}

\section{Applications to the approximate integration}

\begin{lem}\label{Lm_Taylor}
 If $f\in\C$ then
 \[
  \left|\int_Sf(x)dx-\vol(S)f(\barp)\right|\le\frac{1}{2}\df\int_S\|x-\barp\|^2dx.
 \]
\end{lem}
\begin{proof}
 Write Taylor's formula with the remainder depending on the second order
 differential: for any $x\in S$ there exists a point $\xi_x$ belonging to the
 segment with endpoints $x$ and $\barp$ such that
 \begin{equation}\label{Lm_Taylor_pf}
  f(x)-f(\barp)-df(\barp)(x-\barp)=\frac{1}{2}d^2f(\xi_x)(x-\barp,x-\barp).
 \end{equation}
 Since the left--hand side is an integrable function of $x\in S$, then so is
 the right--hand side. Observe that $x\mapsto df(\barp)(x-\barp)$ is an affine
 function, whence by Lemma~\ref{Lm_aff}
 \[
  \int_S df(\barp)(x-\barp)dx=df(\barp)(\barp-\barp)=0.
 \]
 Therefore integrating both sides of~\eqref{Lm_Taylor_pf} and using~\eqref{norm3}
 we get
 \begin{multline*}
  \left|\int_Sf(x)dx-\vol(S)f(\barp)\right|
  =\frac{1}{2}\left|\int_Sd^2f(\xi_x)(x-\barp,x-\barp)dx\right|\\
  \le\frac{1}{2}\int_S\bigl|d^2f(\xi_x)(x-\barp,x-\barp)\bigr|dx
  \le\frac{1}{2}\df\int_S\|x-\barp\|^2dx.
 \end{multline*}
\end{proof}

\begin{thm}\label{Th_T}
 Let $\T$ be a positive linear functional defined on $\C$ such that
 $\T(p)=\int_Sp(x)dx$ for all polynomials $p:\R^n\to\R$ of degree at most 2. Then
 \[
  \left|\int_Sf(x)dx-\T(f)\right|\le\df\int_S\|x-\barp\|^2dx
 \]
 for any $f\in\C$.
\end{thm}
\begin{proof}
 Let $g(x)=\frac{\df}{2}\|x\|^2$, $x\in S$. By Lemma~\ref{Lm_norm} the
 functions  $g+f$ and $g-f$ are convex. Observe that the operator
 $\frac{1}{\vol(S)}\T$ fulfils the assumptions of Theorem~\ref{Th_Had} (for the
 purposes of this proof it is enough to have $\T$ defined on $\C$, the
 assertion of Theorem~\ref{Th_Had} holds also in this setting). Then
 \[
  \vol(S)(g+f)(\barp)\le\T(g+f)\quad\text{and}\quad\vol(S)(g-f)(\barp)\le\T(g-f).
 \]
 Using the first of the above inequalities, Lemma~\ref{Lm_Taylor} and the second
 assertion of Lemma~\ref{Lm_norm} we infer that
 \begin{multline*}
  \int_S(g+f)(x)dx
  \le\vol(S)\bigl(g(\barp)+f(\barp)\bigr)+\frac{1}{2}\dgpf\int_S\|x-\barp\|^2dx\\
  \le\T(g+f)+\df\int_S\|x-\barp\|^2dx.
 \end{multline*}
 Since $\int_Sg(x)dx=\T(g)$, then
 \begin{equation}\label{Th_T_pf}
  \int_Sf(x)dx-\T(f)\le\df\int_S\|x-\barp\|^2dx.
 \end{equation}
 Repeating the above argument for $g-f$ we arrive at
 \[
  -\left[\int_Sf(x)dx-\T(f)\right]\le\df\int_S\|x-\barp\|^2dx,
 \]
 which, together with~\eqref{Th_T_pf}, finishes the proof.
\end{proof}

To approximately compute multiple integrals on simplices using the cubature
rules it is enough to have error bounds of such rules on the unit simplex
$S_1$, i.e. on the simplex with vertices $(0,0,\dots,0)$, $(1,0,\dots,0)$,
$(0,1,\dots,0)$, \dots, $(0,0,\dots,1)$. The integral over any simplex by the
appropriate affine variable interchange can be computed as an integral over
$S_1$. That is why we record below the following

\begin{cor}\label{Cor_S1}
 Let $\T$ be a positive linear functional defined on $\mathcal{C}^2(S_1)$ such
 that $\T(p)=\int_{S_1}p(x)dx$ for all polynomials $p:\R^n\to\R$ of degree
 at most 2. Then
 \[
  \left|\int_{S_1}f(x)dx-\T(f)\right|\le\frac{n^2}{(n+2)!(n+1)}\df
 \]
 for any $f\in\mathcal{C}^2(S_1)$.
\end{cor}
\begin{proof}
 It is well known that $\vol(S_1)=\frac{1}{n!}$. To prove the Corollary by
 Theorem~\ref{Th_T} it is enough to check that
 \[
  \int_{S_1}\|x-\barp\|^2dx=\frac{n^2}{(n+2)!(n+1)},
  \quad\text{where}\;\barp=\left(\frac{1}{n+1},\dots,\frac{1}{n+1}\right).
 \]
 For $x=(x_1,\dots,x_n)\in S_1$ we have
 \begin{multline*}
  \|x-\barp\|^2=\sum_{i=1}^n\left[x_i-\frac{1}{n+1}\right]^2
  =\sum_{i=1}^n\left[x_i^2-\frac{2}{n+1}x_i+\frac{1}{(n+1)^2}\right]\\
  =\sum_{i=1}^n\left[\pi_i^2(x)-\frac{2}{n+1}\pi_i(x)+\frac{1}{(n+1)^2}\right],
 \end{multline*}
 where $\pi_i(x)=x_i$ is the projection to the $i$--th axis.
 Bessenyei~\cite[Lemma~1]{Bes08} computed
 \[
  \int_{S_1}\pi_i(x)dx=\frac{1}{(n+1)!},\quad i=1,\dots,n.
 \]
 Repeating his argument we can easily compute
 \[
  \int_{S_1}\pi_i^2(x)dx=\frac{2}{(n+2)!},\quad i=1,\dots,n.
 \]
 Then
 \begin{align*}
  \int_{S_1}\|x-\barp\|^2dx
  &=\sum_{i=1}^n
     \int_{S_1}\left[\pi_i^2(x)-\frac{2}{n+1}\pi_i(x)+\frac{1}{(n+1)^2}\right]dx\\
  &=\sum_{i=1}^n
     \left[
      \frac{2}{(n+2)!}-\frac{2}{n+1}\cdot\frac{1}{(n+1)!}
       +\frac{1}{(n+1)^2}\cdot\frac{1}{n!}
     \right]\\
  &=n\left[\frac{2}{(n+2)!}-\frac{1}{(n+1)!}\cdot\frac{1}{(n+1)}\right]
   =\frac{n^2}{(n+2)!(n+1)}.
 \end{align*}
\end{proof}

The above Corollary is useful in the approximate integration since an important
class of positive linear functionals which coincide with the integral for
polynomials of degree at most~2 are operators of the form of a~conical
combination of some values of an integrand at appropriately chosen points.
Below we give an example of such cubature operator.

\begin{cor}
 Let $T$ be the unit simplex in $\R^2$,
 \[
  \T(f):=\frac{1}{24}\bigl(f(0,0)+f(0,1)+f(1,0)\bigr)
        +\frac{3}{8}f\left(\frac{1}{3},\frac{1}{3}\right).
 \]
 Then
 \[
  \left|\iint_Tf(x,y)dxdy-\T(f)\right|\le\frac{1}{18}\df
 \]
 for any $f\in\mathcal{C}^2(T)$.
\end{cor}
\begin{proof}
 A cubature operator $\T$ fulfils all the assumptions of
 Corollary~\ref{Cor_S1}, in particular, $\T(p)=\iint_Tp(x,y)dxdy$ for
 any polynomial $p:\R^2\to\R$ of degree at most~2.
\end{proof}

\bibliographystyle{amsplain}
\bibliography{was_own,was_pub}

\end{document}